\documentclass[twoside,a4paper,12pt,centertags,reqno]{amsart} % 13pt displays better than 12pt
\usepackage{amsmath,amssymb,verbatim,vmargin}
\usepackage{color}
\usepackage{tikz}
\usepackage{hyperref}
\hypersetup{colorlinks,bookmarks=true,linktocpage=true,citecolor=blue,linkcolor=magenta}

\usepackage{eulervm}   % formulas round   

\allowdisplaybreaks % align break
\usepackage{color}
\usepackage[all]{xy} % diagrams

\usepackage{mathtools}
\usepackage{MnSymbol} %wedge righthalfcup

\theoremstyle{plain}
\newtheorem{thm}{Theorem}[section]

\theoremstyle{definition}

\newtheorem{rem}[thm]{Remark}

\newtheorem{fact}[thm]{Fact}

\usepackage{enumerate}

\newcommand{\bQ}{{\mathbb Q}}

\newcommand{\bN}{{\mathbb N}}

\def\barint_#1{\mathchoice
            {\mathop{\vrule width 6pt
height 3 pt depth -2.5pt
                    \kern -9.5pt
\intop \kern -4pt}\nolimits_{#1}}%
            {\mathop{\vrule width 5pt height
3 pt depth -2.6pt
                    \kern -6.5pt
\intop \kern -4pt}\nolimits_{#1}}%
            {\mathop{\vrule width 5pt height
3 pt depth -2.6pt
                    \kern -6pt
\intop \kern -4pt}\nolimits_{#1}}%
            {\mathop{\vrule width 5pt height
3 pt depth -2.6pt
          \kern -6pt \intop \kern -4pt}\nolimits_{#1}}}
          
           \def\bariint_#1{\mathchoice
            {\mathop{\vrule width 15pt
height 3 pt depth -2.5pt
                    \kern -15.8pt
\intop \kern -8pt\intop \kern -4pt}\nolimits_{#1}}%
            {\mathop{\vrule width 9pt height
3 pt depth -2.6pt
                    \kern -10.5pt
\intop \kern -8pt\intop \kern -4pt}\nolimits_{#1}}%
            {\mathop{\vrule width 9pt height
3 pt depth -2.6pt
                    \kern -10pt
\intop \kern -8pt\intop \kern -4pt}\nolimits_{#1}}%
            {\mathop{\vrule width 9pt height
3 pt depth -2.6pt
          \kern -8pt \intop \kern -10pt\intop \kern -4pt}
      \nolimits_{  #1}}}

\def\barintlim_#1{\mathchoice
            {\mathop{\vrule width 6pt
height 3 pt depth -2.5pt
                    \kern -8.8pt
\intop \kern -4pt}\limits_{#1}}%
            {\mathop{\vrule width 5pt height
3 pt depth -2.6pt
                    \kern -6.5pt
\intop \kern -4pt}\limits_{#1}}%
            {\mathop{\vrule width 5pt height
3 pt depth -2.6pt
                    \kern -6pt
\intop \kern -4pt}\limits_{#1}}%
            {\mathop{\vrule width 5pt height
3 pt depth -2.6pt
          \kern -6pt \intop \kern -4pt}\limits_{#1}}}
          
           \def\bariintlim_#1{\mathchoice
            {\mathop{\vrule width 15pt
height 3 pt depth -2.5pt
                    \kern -15.8pt
\intop \kern -8pt\intop \kern -4pt}\limits_{#1}}%
            {\mathop{\vrule width 9pt height
3 pt depth -2.6pt
                    \kern -10.5pt
\intop \kern -8pt\intop \kern -4pt}\limits_{#1}}%
            {\mathop{\vrule width 9pt height
3 pt depth -2.6pt
                    \kern -10pt
\intop \kern -8pt\intop \kern -4pt}\limits_{#1}}%
            {\mathop{\vrule width 9pt height
3 pt depth -2.6pt
          \kern -8pt \intop \kern -10pt\intop \kern -4pt}
      \limits_{  #1}}}
          
\renewcommand{\iint}{\int \kern -3pt\int}       

\numberwithin{equation}{section}
\title{Change of variable and discrete Hardy inequality}  

\dedicatory{In memory of Jonathan Sondow (1943-2020)}

\author{Yi C. Huang} 
\address{School of Mathematical Sciences, Nanjing Normal University, Nanjing 210023, People's Republic of China}
\email{Yi.Huang.Analysis@gmail.com}
\urladdr{https://orcid.org/0000-0002-1297-7674}

\date{\today} 

\thanks{Research of the author is partially supported by the National NSF grant of China (no. 11801274).
The author would also like to thank Professor Pascal Lef\`evre (Artois) for helpful communications.}

\begin{document}

\begin{abstract}
For absolutely convergent series we state explicitly a one-sided summation estimate that can be viewed as the discrete analogue of the change of variable formula on the half line.
This estimate is implicit in Pascal Lef\`evre's recent elegant proof of the classical discrete Hardy inequality.
Here we remove a superfluous irrationality condition therein and point out the change of variable character of his approach.
This leads to a simpler, shorter and \textit{bona fide} Ingham type proof of the discrete Hardy inequality, and also provides the optimal constant.
\end{abstract}

\maketitle

\section{Introduction}

Let $p\geq1$. 
Denote by $\ell^p$ the space of sequences of complex numbers $\{a_n\}_{n\geq0}$ such that $\sum_{n\geq0}|a_n|^p<+\infty$,
endowed with the norm 
$$\|\{a_n\}\|_p:=\left(\sum_{n\geq0}|a_n|^p\right)^{1/p}.$$
Recently in \cite{Lef20}\footnotemark
\footnotetext{See also \cite{Lef20UMJ} for a subsequent extension to weighted discrete Hardy's inequalities.}, 
Pascal Lef\`evre gave a short proof of the classical discrete Hardy inequality in the sequence case, that is, the $\ell^p$-norm estimate
\begin{equation} \label{eqn:dhi}
\|\{A_n\}\|_p\leq p'\|\{a_n\}\|_p,\quad p>1.
\end{equation}
Here, $\{A_n\}$ is the Ces\`aro mean of $\{a_n\}$ given by
$$A_n=\frac{1}{n+1}\sum_{k=0}^na_k$$ 
and $p'=\frac{p}{p-1}$ is the optimal constant for \eqref{eqn:dhi}.
This famous and useful inequality can be found for example in Hardy, Littlewood, and P\'olya \cite{HLP67}.
The proof machinery of Lef\`evre is \textit{direct}, in the sense that it does not appeal to the continuous version of Hardy inequality and it also provides the sharp constant.
However, a framework that is able to take care of both the discrete case and the continuous case would be conceptually more desirable.
In this spirit is the supersolution perspective taken by Keller, Pinchover, and Pogorzelski \cite{KelPinPog18} on the discrete Hardy inequality, with improved Hardy weights; 
their result can be regarded as a natural discrete translation of Devyver, Fraas, and Pinchover \cite{DevFraPin14} in the continuum.
These improved discrete inequalities were later reproved by the author in \cite{Hua21}
with motivations from the Machihara-Ozawa-Wadade theory \cite{MacOzaWad17,MacOzaWad19} on Hardy and Rellich inequalities in the framework of equalities.
For an abstract Machihara-Ozawa-Wadade theory and a systematic treatment of Hardy inequalities on homogeneous groups, see Ruzhansky and Suragan \cite{RuzSur19}.
For historical aspects and several other proofs of \eqref{eqn:dhi} and its continuous version, see Kufner, Malgranda, and Persson \cite{KMP06}.
For more recent developments about Hardy inequalities on Euclidean domains with nonempty boundary, 
see for example Avkhadiev \cite{Avk06, Avk13, Avk21}, Goel, Pinchover, and Psaradakis \cite{GoePinPsa22} and the references therein.

Back to Lef\`evre's arguments for \eqref{eqn:dhi}, the following change of variable \textit{estimate} for convergent series is implicit.
We single it out for its potential independent interest.

\begin{fact} \label{fact:cov}
Let $s>0$ and $\|\{a_n\}\|_1<+\infty$. Then
\begin{equation} \label{eqn:cov}
s\sum_{n\geq1}|a_{[ns]}|\leq \sum_{n\geq0}|a_n|.
\end{equation}
Here $[x]$ denotes the largest integer not exceeding $x$.
\end{fact}

More precisely, that \eqref{eqn:cov} holds with $s\in[0,1]\backslash\bQ$ was used in \cite{Lef20} to prove the discrete Hardy inequality \eqref{eqn:dhi}.
In this note we point out that this additional $``\backslash\bQ"$ condition is superfluous. 
Note that \eqref{eqn:cov} also holds trivially when $s\in\{0,1\}$. 

It is plain that \eqref{eqn:cov} is similar to the change of variable \textit{formula} in the continuum
\begin{equation} \label{eqn:covcont}
s\int_0^{+\infty} f(sx)dx=\int_0^{+\infty} f(x)dx, \,\,\forall \,s>0.
\end{equation}
Under the Dirichlet boundary condition $a_0=0$ (see \cite{KelPinPog18} for the interpretation), the equality in \eqref{eqn:cov} is attained when no change of variable happens, i.e., $s=1$.
Moreover, \eqref{eqn:cov} becomes especially clear-cut for integer $s$ and non-increasing $\{|a_n|\}$.

\section{Proofs of Fact \ref{fact:cov}}

It is a crucial observation that it suffices to work with the non-increasing rearrangement of $\{|a_n|\}$.
Denote by $\bN$ the set of positive integers and by $\bN_0=\bN\cup\{0\}$.

\begin{proof}[First proof via counting function (after Lef\`evre)]
For $s>0$ and $m\in\bN_0$, introduce
$$I_m(s)=\{n\in\bN:\,\,[ns]=m\}=[m/s,(m+1)/s)\cap\bN$$
and
$$J_m(s)=[0,m/s)\cap\bN.$$
Denote by $\sharp$ the counting function. Now, following the arguments in \cite{Lef20}, 
$$\begin{aligned}
\sum_{n\geq1}|a_{[ns]}|&=\sum_{m\geq0}\sum_{n\in I_m(s)}|a_{[ns]}|\\
&=\sum_{m\geq0}|a_m|\sharp I_m(s)\\
&=\sum_{m\geq0}|a_m|\left(\sharp J_{m+1}(s)-\sharp J_m(s)\right).
\end{aligned}$$
Using Abel transform, $\sharp J_m(s)< m/s$ for $m\geq1$ and monotonicity of $\{|a_n|\}$, one has
$$\begin{aligned}
\sum_{m\geq0}|a_m|\left(\sharp J_{m+1}(s)-\sharp J_m(s)\right)&=\sum_{m\geq0}\sharp J_{m+1}(s)(|a_m|-|a_{m+1}|)\\
&< s^{-1}\sum_{m\geq0}(m+1)(|a_m|-|a_{m+1}|).
\end{aligned}$$
The estimate \eqref{eqn:cov} is thus proved after another use of Abel transform.
\end{proof}

\begin{rem}
In this proof we claim by no means any originality over Lef\`evre's arguments.
What we do is to \textit{not} compute $\sharp I_m(s)$, and just leave it implicit.
\end{rem}

\begin{proof}[Second proof via change of variable formula]
As observed above, we can work with decreasing $\{|a_n|\}$.
Then, \eqref{eqn:cov} simply follows from \eqref{eqn:covcont} applied to $f(x)=|a_{[x]}|$.
\end{proof}

\section{Change of variable and discrete Hardy inequality}

To end our discussions, let us recall Lef\`evre's elegant proof of \eqref{eqn:dhi}. 
By an Ingham type representation $A_n=\int_0^1a_{[(n+1)s]}ds$
and triangle inequality, one has
\begin{equation} \label{eqn:Ingham}
\left(\sum_{n=0}^N|A_n|^p\right)^{1/p}\leq \int_0^1\left(\sum_{n\geq1}|a_{[ns]}|^p\right)^{1/p}ds.
\end{equation}
Inserting the \textit{pointwise} inequality \eqref{eqn:cov} for $\{|a_{[ns]}|^p\}$ into \eqref{eqn:Ingham}, computing the resulted integral and letting $N\rightarrow+\infty$ 
then complete the proof of \eqref{eqn:dhi}.

In conclusion, although Lef\`evre wrote in Remark (2) of \cite{Lef20} that ``$\cdots$ there is no change of variable in our argument",
our interpretation in this note somehow clarifies that the ingredient of his approach is indeed a change of variable in the discrete setting.
Furthermore, our proof of the discrete Hardy inequality \eqref{eqn:dhi} via the second approach to \eqref{eqn:cov}  is simpler and shorter than the one given in \cite{Lef20}.

\bigskip

\section*{\textbf{Compliance with ethical standards}}

\bigskip

\textbf{Conflict of interest} The author declares that there is no conflict of interest.

\bigskip

\textbf{Availability of data and material} Not applicable.

\bigskip

\bibliographystyle{alpha}

\bibliography{Hua-DiscHardyCount} 
 
\end{document}